\vspace{10pt}
\documentclass{llncs}
\usepackage{amssymb}

\newtheorem{thm}{Theorem}

\def\x{$\hfill\rlap{$\sqcup$}\sqcap$\bigskip}

\usepackage[fleqn]{amsmath}

\newtheorem{cor}[thm]{Corollary}

\newcommand{\F}{{\mathbb{F}}}

\begin{document}

\title{Quadratic Binomial APN Functions and Absolutely Irreducible
Polynomials }
\author{Eimear Byrne and Gary McGuire\thanks{Research 
supported by the Claude
Shannon Institute, Science
Foundation Ireland Grant 06/MI/006.
Email {\tt gary.mcguire@ucd.ie}}}

\institute{School of Mathematical Science, University College Dublin, Ireland}


\maketitle

\begin{abstract}
\noindent 
We show that many quadratic binomial  functions of the form
$cx^{2^i+2^j}+dx^{2^u+2^v}$ ($c,d\in GF(2^m)$) are not APN infinitely often.
This is of interest in the light of recent discoveries of new families of quadratic binomial
APN functions.
 The proof uses the Weil bound from algebraic geometry.
\end{abstract}


\section{Introduction}
Let $K:=GF(q)$ be the finite field with $q$ elements.
Let $n$ be a positive integer, and let $L$ be an extension of $K$ of degree $n$.
We consider polynomials with coefficients in $K$, as polynomial functions on $L$. Let $f \in K[x]$.
The function $f: L\longrightarrow L$ is called {\em perfect nonlinear} (PN) on $L$ if for every $a \in L^*$, $b \in L$, there is at most one solution $x \in L$ to the equation
\begin{equation}\label{eqapn}
f(x+a) - f(x) = b.
\end{equation}
There are no perfect nonlinear functions on fields of characteristic 2, since whenever $x$ is a solution to (\ref{eqapn}), then so is $x+a$. If for every $a \in L^*$, $b \in L$, there are at most two solutions $x \in L$ to (\ref{eqapn}), then we say that $f$ is {\em almost perfect nonlinear} (APN) on $L$.
Due to the connections to coding theory and cryptography, APN functions over fields of characteristic 2
are more widely studied. For the remainder we assume that $q=2^m$ for some positive integer $m$.

An equivalent definition is to say that a function $f$ is APN on $L$
if the set
$$D_a = \{ f(x)+f(x+a) : x \in L\}$$
is as large as possible (namely $q^{n}/2$) for every nonzero $a \in L$. 
$D_a$ is called the {\em differential} of $f$ at $a$. 
By definition, an APN function provides best possible resistance to a {\em differential attack} when used as an S-box of a block cipher, since then 
the plaintext difference $a=x+y$ yields the ciphertext difference $b=f(x)+f(y)$ with least probability.

Until 2006, the list of known affine inequivalent APN functions on $K=GF(2^m)$ was rather short:

\begin{center}
\begin{tabular}{|l|l|l|l|}
\hline
 $x^d$ & exponent $d$ & constraints & \\
 \hline
 Gold & $2^r+1$ & $(r,m) = 1$ & \cite{G,N} \\
 \hline
 Kasami-Welch & $2^{2r}-2^r+1$ & $(r,m)=1,m$ odd & \cite{K} \\   
 \hline
 Welch &$2^r+3$ & $n=2r+1$  & \cite{Dobb}\\
 \hline
 Niho & $2^r+2^{r/2} - 1$      & $ m = 2r+1, t$ even & \cite{Dobb1} \\
                    & $2^r+2^{(3r+1)/2}-1$  &$m = 2r+1, t$ odd &  \\
\hline
Inverse      &$2^{2r}-1$& $m=2r+1$ & \cite{BD,N} \\
\hline
Dobbertin &$2^{4r}+2^{3r}+2^{2r}+2^r-1$& $ m = 5r$ & \cite{Dobb2}\\
\hline
\end{tabular} 
\end{center}

It was conjectured that the list was complete, up to equivalence. 
Motivated by this question, several authors have considered when
a linear combination of Gold functions  could be APN.
In \cite{BCC+} the authors show that a polynomial of the form
$$f(x) = \sum_{j \in J} f_{j} x^{2^j+1} \in K[x], $$
$J \subset {\mathbb N}$, with at least two nonzero $f_j$ is not APN on $K$ (and therefore not on any extension $L$ of $K$).
The main result relies on proving that the polynomial
$$\frac{f(x)}{x^2} = \sum_{j \in J} f_{j} x^{2^j-1}$$
is a permutation polynomial on K. While the authors in \cite{BCC+} use Hermite's criterion to establish the permutation property, in fact this was proved sometime before in \cite{P} by different techniques.
In any case, this method cannot be extended to more general quadratics, with exponents of the form
$2^j+2^k$ instead of $2^j+1$.
It is this more general form of exponent that we will consider in this paper. 
We will study when a linear combination of two Gold functions with these more general
exponents can be APN.

This case with more general quadratic exponents is more interesting, because such functions \emph{can} indeed be APN.  
In \cite{EKP}, the first example of an APN function
not equivalent to any of the above list appeared. The function
 \begin{equation}\label{EKP}
            x^3 + u x^{36} \in GF(2^{10})[x],$$
\end{equation}
where $u \in \omega GF(2^5)^*\cup \omega^2 GF(2^5)^* $ and 
$\omega$ has order 3 in $GF(2^{10})$ is APN on $GF(2^{10})$. 
This function has the additional property of being {\em crooked}. The function $f$ is called
crooked if $D_a$ forms an affine hyperplane over $GF(2)$. This will certainly be the case
for any APN function $f$ for which $f(x)+f(y)+f(x+y)$ is $GF(2)$-bilinear, in which case $f$ is called a {\em quadratic} APN function. Crooked functions have connections with other combinatorial objects,
such as distance-regular graphs \cite{dCMM,dCvD}. It is now known that all crooked monomials or binomials are quadratic \cite{BK,K05}.
  
Since the emergence of this first sporadic example, there are now several more known infinite families of inequivalent APN functions.

\begin{center}
\begin{tabular}{|l|l|l|l|}
\hline
& $ f(x)$ & constraints &  \\
\hline
   (1) & $x^{2^i+1} + u x^{2^{k+i}+2^{k(r-1)}}$ & $ m=rk \in \{ 3k, 4k\}, (r,s)=(s,k)=(r,k)=1$, & \cite{BCFL,BCL1}\\
     
                                                                  &        &   $r|(k+s),$ and $u$ is a $2^k-1$-th power in $K$ &\\
\hline                                                                       
 (2) & $ux^{2^{-k}+2^{k+s}}+u^{2^k}x^{2^{s}+1}+vx^{2^{k+s}+2^s} $ &  $ m=3k, (3,s)=(s,k)=(3,k)=1$, $uv \neq 1$
 & \cite{BB1+}\\
     
                                                                     &     &   $3|(k+s),$ and $u$ is a $2^k-1$-th power &\\
\hline
(3) & $ bx^{2^s+1}+b^{2^k}x^{2^{k+s}+2^k}+cx^{2^{k}+1} $ & $m=2k$,  $r_j \in GF(2^k)$, & \cite{BB1+}\\
 & $+ \sum_{j=1}^{k-1}r_j x^{2^{j+k}+2^i}$                           & $b,c \notin GF(2^k)$&\\    
\hline
(4) & $x^3 + {\rm Tr}(x^9)$ & & \cite{BCL2} \\                                                                                                    
\hline
(5) & $ u^{2^k} x^{2^{-k}+2^{k+s}} + u x^{2^s+1} + v x^{2^{k+s}+2^s}$ & $m=3k, (3,s)=(s,k)=(3,k)=1$, $3|(k+s),$ & \cite{BB2+}\\
                                                                       &   &  $u$ primitive and $v \in GF(2^k)$ &\\
\hline
(6) & $ u^{2^k} x^{2^{-k}+2^{k+s}} + u x^{2^s+1} + v x^{2^{-k}+1} $ & $m=3k, (3,s)=(s,k)=(3,k)=1$, $3|(k+s),$ 
& \cite{BB2+}\\
  & $+ wu^{2^k+1}x^{2^{k+s}+2^s}$                                                  &  $u$ primitive and $w,v \in GF(2^k)$, $wv \neq 1$ &\\
\hline                                                                          
\end{tabular}\label{tab2}
\end{center} 

Observe that the new families are all quadratic APN functions, and
have exponents of the more general form $2^j+2^k$.
This is partially explained by the fact that proving the APN property
for quadratics seems to be easier than for arbitrary polynomials.
Further, at the time of writing,
not a single example of a non-monomial non-quadratic APN function is known.

Family (1) generalizes a second example given in \cite{EKP}, for the case $m=12$. Family (2) contains Family (1) as a subclass for the case $m=3k$ when $v=0$. Family (6) contains Family (5) as a subclass when $w=0$.
Most of these families include some of the inequivalent polynomials listed in \cite{Dlist}. 
In fact Family (2) (with $r_i=0$) contains a class equivalent to the trinomial family given in \cite{BC}.
In the same paper, based on a construction presented in \cite{Dlist}, the authors present a family of hexanomials on $GF(2^{2k})$:
$$x(x^{2^i}+x^{2^k}+cx^{2^{i+k}})+x^{2^i}(c^{2^k}x^{2^k}+bx^{2^{i+k}})+x^{2^{i+k}+2^k} $$
with $k\geq 3, (i,k)=1, b \notin GF(2^k)$, and show that such a hexanomial is APN as long as
$$p(x)=x^{2^i+1}+cx^{2^i}+c^{2^k}x +1 $$
is irreducible over $GF(2^{2k})$. It is not clear that the polynomial $p(x)$ can always be shown to be irreducible over $GF(2^{2k})$ for some choice of $c$. For this reason we do not yet classify it as an infinite family of APN functions. However, the authors of \cite{BC} have checked by computer that for $6<m<26$, several such $c$ exist (about 3/10 of all field elements).

An important component of the study of APN functions is the notion of equivalence of APN functions.
The most prominent of these relations are {\em extended affine equivalence} (EA) and CCZ-equivalence 
\cite{CCZ98}. A pair of functions are EA equivalent if one can be obtained from the other by affine permutations;
more precisely, $f \equiv_{EA} g$ if there exist affine permutations $A_1,A_2$ and an affine map $A$ satisfying
$f=A_1 \circ g \circ A_2 + A$. If $f \equiv_{CCZ} g$ then the graph of $f$ can be obtained from the graph of $g$ by an affine permutation. CCZ-equivalence generalizes affine equivalence and can be expressed in terms of coding theory. In fact two functions are CCZ -equivalent if and only if some corresponding linear codes are equivalent (cf \cite{BB1+,Dlist}). CCZ-equivalent functions have the same {\em differential uniformity} (which is 2 for APN functions) and the same {\em nonlinearity} and hence offer the same resistance to linear and differential attacks. In general, proving CCZ-inequivalence is very difficult and in several instances we rely on computing power to return inequivalent functions from different families. Although typical computations in establishing CCZ-equivalence become intractable even for relatively small values of $m$, a recent approach (see \cite{BB2+,CNS} for further details)
has shown that Families (6) and (7) are CCZ-inequivalent to any previously known families (working over $GF(2^{12})$).  

Recall that $K=GF(2^m)$.
For each family listed above, there is a functional (or equational) relationship between the parameter $m$ and the parameters $k,s$ etc. appearing in the function.  
In other words, the form of the function depends on $m$.
This includes the $x^3+Tr(x^9)$ function, since the trace term depends on $m$.
This is not the case for the Gold and Kasami-Welch functions; any fixed Gold or Kasami-Welch is 
APN on infinitely many extensions of $K=GF(2)$.
We will say that a function $f$ that is defined over $K$  is {\em APN infinitely often}
if $f$ is APN on $K$ and APN on an infinite number of extensions of $K$.
One way to tackle the classification problem of APN functions is to determine classes of APN functions that are not APN infinitely often. 
This problem has been studied for monomial functions in 
\cite{JMW} and \cite{Jed}, and for arbitrary polynomials more recently by Rodier \cite{R}.
We  take this approach in this paper, focussing on  quadratic binomials. All approaches to this invoke the Weil bound and its generalizations. 

We will show that many classes of quadratic binomials defined on $K$ are 
not APN infinitely often.  
A summary of our results is in Corollary \ref{weilagain}.
One consequence of our results is  that the APN binomials of Family (1) are not APN infinitely often, for $i > 1$. Using Frobenius arguments and singularities of a curve, we also show that the sporadic APN polynomial $x^3+ux^{36} \in GF(2^{10})[x]$ is not APN infinitely often. 
  
We conjecture that the Gold and Kasami-Welch monomial functions are the only APN
functions that are APN infinitely often.


\section{The Weil Bound}\label{prelim}
A theorem of Weil gives an upper bound on the number of rational points
of an absolutely irreducible curve \cite{W48}. In fact there have been a number of improvements to this
bound since (see \cite{FJ,LN}). However, for the purposes of this paper, we only require the fact that the number of rational points over $L=GF(q^n)$ of an absolutely irreducible curve exceeds its degree for sufficiently large $n$. 

\begin{theorem} \cite{LN}
Let $A(x,y)$ be an absolutely irreducible polynomial of degree $d$ with 
coefficients in $K=GF(q)$. Then the number of points $N_{n}$
on the affine curve $A(x,y)=0$ over $L=GF(q^n)$ satisfies
$$|N_{n} -(q^{n}+1)| \leq 2(d-1)(d-2)\sqrt{q^n} +C_d,$$ 
where $d$ denotes the degree of $A(x,y)$ and $C_d$ is a constant that depends only on $d$.
\end{theorem}

\begin{cor}\label{weil}
Let $G(x,y) \in K[x,y]$ have an absolutely irreducible factor in $K[x,y]$.
Then $G(x,y)$ has rational points over $L$ off the line $x=y$ for all
$n$ sufficiently large.  
\end{cor}

\noindent
Proof:
Let $A(x,y)$ be an absolutely irreducible factor of
$G(x,y)$ in $K[x,y]$ and suppose that
$A(x,y)$ has degree $d$. Then there are at most
$d$ points $(x,y)$ on $A(x,y)$ with $x=y$.
From the Weil bound, for $n$ sufficiently large, the total number of rational
points of the affine curve $A(x,y)=0$ exceeds $d$, which is the upper bound of the number of rational
points of $A(x,y)=0$ on the line $x=y$.
\qed
\bigskip

We apply this result as follows. Given a function $f$ on $L$, define
$$ \Delta_f(x,y) := f(x+y)+f(x)+f(y).$$
Suppose that $f$ is quadratic.
Then clearly $\Delta(a,a)= \Delta(0,a)=0$
for all $a \in L$, so $f$ is APN on $L$ if and only if  $\Delta(x,y)$ has no
rational points over $L$ off the line $x=y$.
Therefore, in light of Corollary \ref{weil}, 
we will study $f$ for which $\Delta_f(x,y)$ has an absolutely irreducible factor.
 Such $f$ will not be APN infinitely often.

\section{Quadratic Binomials}\label{QB}

We will show that many binomial quadratic functions of the form
$cx^{2^i+2^j}+dx^{2^u+2^v}$ ($c,d\in K$) are APN on at most a finite number
of extensions of $K$.

\medskip

\noindent
Observe first that every quadratic binomial in $K[x]$ is affine equivalent to a
function $f$ of the form
\begin{eqnarray*}
   f(x)=x^{2^i+1}+\delta x^{2^{s}(2^{t}+1)}
\end{eqnarray*}
for some nonzero $\delta$ in $K$, where $i\geq 1$, $t\geq 1$, $s\geq 0$.
Then
\begin{eqnarray*}
   \Delta_f(x,y) &: = & f(x+y)+f(x)+f(y) \\
                 & = & x^{2^i}y+xy^{2^i}+ \delta
                 (x^{2^t}y+x y^{2^t})^{2^{s}}.
\end{eqnarray*}
When $s=0$, for any $x \in L^*$, $\Delta_f(x,a) =0 $ if and only if
$$\frac{\Delta_f(x,a) }{xa} = x^{2^i-1}+\delta x^{2^t-1} + a^{2^i-1}+\delta a^{2^t-1} =0.$$
This equation has only the solution $x=a$ in $L$ iff $x^{2^i-1}+\delta x^{2^t-1} $ is a permutation polynomial on $L$, which is  false by the result of \cite{BCC+} mentioned in the introduction.
We deduce that $f$ is not APN on $K$ (or any extension) if $s=0$. 
We therefore assume for the remainder that $s >0 $.
\bigskip

\begin{thm}\label{corweil}
Define $F(x,y)$ by
\begin{eqnarray}\label{eqFq}
    F(x,y)&:=&\frac{\Delta_f(x,y)}{xy} =  x^{2^i-1}+y^{2^i-1}+
   \delta (xy)^{2^s-1}
   (x^{2^t-1}+y^{2^t-1})^{2^{s}}.
\end{eqnarray}
 If $F(x,y)$ has an absolutely irreducible factor over $K$, then  $f(x)$ is not APN infinitely often.
\end{thm}

Proof:
 If $F(x,y)$ has an absolutely irreducible factor over $K$, then  $f(x)$ is not APN on $L$ for all
$n=[L:K]$ large enough by Corollary \ref{weil} and the subsequent remarks.
\x

We next observe some obvious factors of  $F(x,y)$.

\begin{lemma}\label{div}
Let $d=gcd(i,t)$.
Then $x^{2^d-1}+y^{2^d-1}$ divides $F(x,y)$.
\end{lemma}

Proof:  Let $\alpha\in \F_{2^d}$ be nonzero. Observe that
$F(x,\alpha x)=0$. This means that $y-\alpha x$ divides $F(x,y)$,
and so
\[
   \prod_{\alpha \in \mathbb{F}_{2^d}, \alpha\not= 0} (y-\alpha x) =
   y^{2^d-1}+x^{2^d-1}
\]
divides $F(x,y)$ in $K[x,y]$. \x

 Let  $U(x,y):=x^{2^d-1}+y^{2^d-1}$ be this obvious factor, and define
 \[
 H(x,y):=\frac{F(x,y)}{U(x,y)}.
 \]

   Then
\begin{eqnarray*}
   H(x,y)&= & \frac{x^{2^i-1}+y^{2^i-1}}{x^{2^d-1}+y^{2^d-1}}+
   \delta (xy)^{2^{s}-1}\frac{x^{2^t-1}+y^{2^t-1}}{x^{2^d-1}+y^{2^d-1}}(x^{2^t-1}+
y^{2^t-1})^{2^s-1}\\
         & = & \prod_{\alpha \in {\cal I}}(x+\alpha y) +
        \delta (xy)^{2^{s}-1}\prod_{\beta \in {\cal B}}(x+\beta y)
        \prod_{\gamma \in {GF(2^t)}}(x+\gamma y)^{2^{s}-1}\\
         & = & H_{2^i-2^d} + H_{2^{t+s}+2^{s}-2^d-1}, \:\:\:\:\:\:\:\:\:\:\:\:\:\
         \:\:\:\:\:\:\:\:\:\:\:\:\:\:\:\:\:\:\:\:\:\:\:\:\:\:\:\:\:\:\:\:\:\:\:\:\:\:\:\:\:\:\:\:\:\:\:\:\:\:\:\:\:\:\:\:\:\:\:\:\:\:\:(4)
\end{eqnarray*}
with ${\cal I} = GF(2^i) \backslash GF(2^d)$ and  ${\cal B} =
GF(2^t)\backslash GF(2^d)$. The polynomials $ H_{2^i-2^d}$, $H_{2^{t+s}+2^{s}-2^d-1} \in K[x,y]$
are homogeneous of degrees $2^i-2^d$ and $2^{t+s}+2^s-2^d-1$ respectively.
\bigskip

We conjecture the following:
\vspace{3mm}

\noindent{\bf Conjecture 1}\;
     The polynomial $H(x,y)$
     is absolutely irreducible.
\vspace{3mm}

Any $f(x)$ for which Conjecture 1 holds is not APN infinitely often, by Theorem \ref{corweil}.

We have strong  results in support of our conjecture. 
We will show that if $i$ does not divide $t$, which is equivalent to the two  homogeneous polynomials $H_{2^i-2^d}$ and $H_{2^{t+s}+2^{s}-2^d-1}$ (that $H$ is the sum of) being non-constant,  then $H$ is  absolutely irreducible. Since $s>0$, if $H$ is not a sum of two non-constant homogeneous polynomials then $H_{2^i-2^d}=1$ and $d=i$. We are not able to prove this second case completely, however we derive some constraints on $i,t,s$ from which we can deduce that $H$ is absolutely irreducible in many cases.

\subsection{The Case $i$ Does Not Divide $t$}

We first prove a general lemma.

\begin{lemma}\label{lemabir}
Let $k$ be a field.  Let $G(x,y)\in k[x,y]$ be the sum of
two nonconstant homogeneous polynomials, i.e.,  $G=G_a+G_b$ where $G_i$ is
homogeneous of degree $i$, and $1\leq a<b$.
Suppose $(G_a,G_b)=1$ and either $G_a$ or $G_b$ factors into
distinct linear factors over $\overline{k}$.
Then $G$ is irreducible over $k$.
\end{lemma}

\noindent Proof:
Suppose not, say $G=WV$.  Write $W$ and $V$ as a sum of homogeneous parts, say
$$
G=WV=(W_e+W_{e+1}+\cdots +W_r)(V_f+V_{f+1}+\cdots +V_s).
$$
Note that $G$ cannot have a homogeneous factor as $(G_a,G_b)=1$, so we have
$W_e\not= 0$, $W_r\not=0$, $V_s\not=0$, $V_f\not=0$, $e<r$, $f<s$.
Then $a=e+f$, $d=r+s$, and $a\leq d-2$.

Suppose first that $G_b$ has distinct linear factors, which implies $(W_r,V_s)=1$.
The term of degree $b-1$ in $G$ is $0=W_r V_{s-1}+W_{r-1}V_s$,
so $W_r$ divides $W_{r-1}V_s$.
As $(W_r,V_s)=1$ we have $W_r|W_{r-1}$ which is impossible by degree
considerations unless $W_{r-1}=0$.
Similarly $V_{s-1}=0$.
Applying this argument successively to terms of degree $b-2$, $b-3$, ...,
we obtain $W_j=0$ for all $j<r$ (and $V_j=0$ for $j<s$),
which means that $W$ is homogeneous, a contradiction.

Secondly we consider the case that $F_a$ has distinct linear factors,
which implies $(W_e,V_f)=1$.
The degree $a+1$ term in $F$ is $0=W_eV_{f+1}+W_{e+1}V_f$, so $W_e$ divides $W_{e+1}V_f$
(and $V_f$ divides $W_eV_{f+1}$).
Since $(W_e,V_f)=1$ we get $W_e|W_{e+1}$ (and $V_f|V_{f+1}$).
Next (assuming $b>a+2$)
the degree $a+2$ term is $0=W_eV_{f+2}+W_{e+1}V_{f+1}+W_{e+2}V_f$,
so $W_e|W_{e+2}$.
Continuing in this way we obtain $W_e|W_j$ (and $V_f|V_j$) for all $j$.
But then $W_e$ is a homogeneous factor of $F$, a contradiction.
\x

\begin{thm}\label{quadf}
   Let $F(x,y) \in K[x,y]$ be defined as in Equation $(\ref{eqFq})$.
   If $i$ does not divide $t$ then 
   Conjecture 1 holds.
\end{thm}

\noindent Proof:
Let $d = (i,t)$. Then $F(x,y) = (x^{2^d-1}+y^{2^d-1})H(x,y)$, where
$H(x,y)= H_{2^i-2^d} + H_{2^{t+s}+2^{s}-2^d-1}$, is defined as in (2). 
If neither $H_{2^i-2^d}$ or $H_{2^{t+s}+2^{s}-2^d-1}$ is constant,
then they are relatively prime and
$H_{2^i-2^d}$ is a product of distinct linear
factors.  We may apply Lemma \ref{lemabir} and conclude
that $H$ is absolutely irreducible.
Our assumption that $s>0$ means that $H_{2^{t+s}+2^{s}-2^d-1}$ is not constant 
(which occurs if and only if $d=t$ and $s=0$).
Finally we note that $H_{2^i-2^d}$ is constant if and only if $d=i$.

\x

\subsection{The Case $i$ Divides $t$}

For $i$ a divisor of $t$, $H(x,y)$ has the form $1+H_{2^{s+t}+2^s-2^i-1}(x,y)$.
We apply singularity analysis and a little Galois theory to establish absolute irreducibility.
It is straightforward to show that the affine curve  $H(x,y)$ is non-singular. 
We consider the homogenized 
projective curves $H(x,y,z)=0$, $U(x,y,z)=0$ and $F(x,y,z)=0$ where
\begin{eqnarray*}\label{defH2}
     H(x,y,z)=z^{2^{t+s}+2^s-2^i-1} + H_{2^{t+s}+2^s-2^i-1}(x,y),
\end{eqnarray*}     
$U(x,y,z)=U(x,y)$ as $U(x,y)$ is already homogeneous, and $$F(x,y,z)=U(x,y,z)H(x,y,z).$$

Let $m_P(F)$ denote the multiplicity of the point $P$ on the curve $F$, etc.

\begin{lemma}\label{mults}
Continue the above notation.  The points in $\mathbb{P}^2(\bar K)$ on $H(x,y,z)$
are of the types given in the following table.
\end{lemma}

\begin{center}
\begin{tabular}{|c|c|c|}
\hline
   $P$  &   & multiplicity \\
\hline
 $[a:b:1]$   &  $H(a,b,1)=0$    & $1$       \\
 $[1:b:0]$   &  $b \in GF(q^i)$ & $2^s-1$       \\
 $[1:b:0]$   &  $b \in GF(q^t) \backslash GF(q^i)$ & $2^s$\\

 $[0:1:0]$   &                  & $2^s-1$       \\
 \hline
\end{tabular}
\end{center}
\bigskip

Proof:
Let $P=[1:b:0]$, $b \in \bar K$. Then 
\begin{eqnarray*}
F(x+1,y+b,z) & = & U(x+1,y+b,z)H(x+1,y+b,z)\\
                    & =  & F_0 +F_1 + F_2 \cdots \\
                    & =  & H_0U_0 + (H_0U_1+H_1U_0) + \cdots  
\end{eqnarray*}
and on the other hand
\begin{eqnarray*}
    F(x+1,y+b,z) &=& z^{2^{t+s}+2^s-2^i-1}\left( (x+1)^{2^i-1}+(y+b)^{2^i-1}\right)\\
    &+&
   \delta (x+1)^{2^s-1}(y+b)^{2^s-1}
   \left((x+1)^{2^t-1}+(y+b)^{2^t-1}\right)^{2^{s}}.
\end{eqnarray*}

Now equating coefficients gives
\begin{eqnarray*}
F_0   =  H_0U_0 =  F(1,b,0)  &=&     
\delta b^{2^s-1}
   (1+b^{2^t-1})^{2^{s}}
\end{eqnarray*}
and
\begin{eqnarray*}
F_1  & = &   H_0U_1 + H_1U_0 \\ 
         & = &   x( \delta b^{2^s-1}
   (1+b^{2^t-1})^{2^{s}}) + y( \delta b^{2^s-2}
   (1+b^{2^t-1})^{2^{s}}) .
\end{eqnarray*}

Suppose that $b \neq 0$. Then we may write:
\begin{eqnarray*}
F_1 &=& (x+b^{-1}y)F_0,\\
F_2 &=& (x^2+b^{-1}xy +b^{-2}y^2)F_0,\\
\vdots & = & \vdots \\
F_{2^s-1} & = & (x^{2^s-1}+b^{-1} x^{2^s-2}y + \cdots b^{-(2^s-1)}y^{2^s-1})F_0,\\
F_{2^s} & = & \delta b^{2^s-1} (x + b^{2^t-2}y)^{2^s}.
\end{eqnarray*} 

It follows that, if $P$ is a point on $F$, then $P$ is singular and $m_P(F) = 2^s$. 
Now $P$ is a point of $U(x,y,z)$ if and only if $b \in GF(2^i)^*$, in which case
$m_P(U) = 1$ and $m_P(H) = 2^s-1$; otherwise $m_P(H)=2^s$.
Finally, $[1:0:0]$ and $[0:1:0]$ are points of $H$ of multiplicity $2^s-1$.
\x

We next make two observations on reducibility using the Frobenius automorphism.
We will combine these with our singularity analysis to determine further conditions
that guarantee the absolute irreducibility of $H$.

\begin{lemma}\label{rdivn}
Let
$h(x,y,z)$ be an irreducible homogeneous polynomial of degree $d$ over $K=GF(q)$, 
with leading coefficient $1$ with respect to some monomial order.
If $h$ has an absolutely irreducible factor with leading coefficient $1$
defined over $L=GF(q^n)$
(and no proper subfield) then $n$ divides $d$.
\end{lemma}

Proof.  
Let $h_1$ be an absolutely irreducible factor of $h$ 
with leading coefficient $1$ defined over $L$
(and no proper subfield) of degree $d_1$.
It follows that no Galois conjugate of $h_1$ is a scalar multiple of $h_1$.
The product of the Galois conjugates of $h_1$, the polynomial 
$$\hat{h} = \prod_{\sigma \in Gal(L:K)} \sigma(h_1),$$
has degree $nd_1$.
Moreover, $\hat{h}$ has coefficients in $K$, so $\hat{h}$ divides $h$ in $K[x,y,z]$.
Since $h$ is irreducible in $K[x,y,z]$, $h=\hat{h}$,
so then $d=nd_1$ and the result follows.
\x

\begin{lemma}\label{mults}
Let $h(x,y,z)$ be an irreducible homogeneous polynomial of degree $d$ in over $K=GF(q)$, with leading
coefficient $1$ with respect to some monomial order.
Let $P$ be a point on $h$ of multiplicity $m_0$, and
suppose the coordinates of $P$ lie in an extension of $K$ of degree $r$.
If $h$ factors into $n$ absolutely irreducible factors  over
$L=GF(q^n)$ where $gcd(n,r)=1$, then $n$ divides $m_0$.
\end{lemma}

Proof.
As in Lemma \ref{rdivn}, 
let $h=h_1 \ldots h_n$ be the factorization of $h$ into
absolutely irreducible factors over $L$,
where each $h_i$ has degree $d/n$.

Since $(r,n)=1$, from the Chinese Remainder Theorem there exists
an integer $a$ such that $a\equiv 0 \pmod r$ and $a\equiv 1 \pmod n$.
Let $\sigma$ be the automorphism in the Galois group of
$GF(q^{nr})$ over $K$
given by $\sigma (z) = z^{q^a}$. Then $\sigma$ fixes $GF(q^{r})$ element-wise and 
is a generator for the Galois group of $L$ over $K$.
Then $\sigma$ fixes $P$ and acts transitively on the $h_i$.
It follows that $P$ has the same multiplicity on each $h_i$,
and so $m_0=n \cdot m_P(h_1)$.
\x

\begin{theorem}
    Let $i,t,s$ be positive integers with $i |t$. Let $H(x,y)$ be defined as in $(4)$. 
    Suppose that $H(x,y)$ is 
    irreducible over $K=\F_q$. If either $(t,2^s-1 )=1$ or $(t-i,s) = 1$ then 
    Conjecture $1$ holds.
\end{theorem} 

Proof.
Of course $H(x,y)$ is irreducible if and only if the homogenization $H(x,y,z)$ is irreducible.
Suppose that $H(x,y,z)$ has a proper factorization into absolutely irreducible factors over $L=GF(q^n)$. Now $H$ has the singular point $[1:0:0]$ over $K$ of multiplicity $2^s-1$, 
so from Lemma \ref{mults} $n | (2^s-1)$.   We must show that $n=1$.

$H$ has also has singular points $[1:a:0]$ with $a \in GF(q^t)$ of multiplicity $2^s$. Moreover, if $(t,2^s-1)=1$ then $(t,n)=1$ as $n$ divides $2^s-1$.
From Lemma \ref{mults} $n | 2^s$, which together with 
$n | (2^s-1)$ forces  $n=1$ and so $H$ is absolutely irreducible.

For the last part, observe that as in the proofs of Lemma \ref{rdivn} and Lemma \ref{mults}, $n$ divides $2^{s+t}+2^s-2^i-1$, the degree of $H$.
Then $n$ divides $2^{s+t}-2^i = 2^i(2^{s+t-i}-1)$, and since $n$ is odd, $n | (2^{s+t-i}-1)$. If $(t-i,s)=1$, this again forces $n=1$.
\x

We now summarize all our results.


\begin{cor}\label{weilagain}
Suppose that
$f = x^{2^i+1}+\delta x^{2^s(2^t+1)}$ where $\delta \in K$, 
$i,t,s >0$. Let $H(x,y)$ be defined as in $(4)$. Suppose that any of the following 
conditions holds:
\begin{enumerate}
    \item
        $i$ does not divide $t$,
    \item   
       $(t,2^s-1)=1$ and $H(x,y)$ is irreducible over $K$,
    \item 
      $(t-i,s)=1$  and $H(x,y)$ is irreducible over $K$.  
 \end{enumerate}       
 Then $f$ is not APN infinitely often.
\end{cor}
 
 \bigskip
 
 Finally, we present some applications of our results.
 
 \bigskip
 
\noindent {\bf Example 1.} We now prove that the binomials of
Family (1) in the introduction are not APN infinitely often.
For $r=3,4$, the binomials of \cite{BCFL,BCL1} in Family (1), namely,
$$x^{2^i+1} + ux^{2^{k+i}+2^{(r-1)k}},$$
defined on $K=GF(2^{rk})$ 
are APN on $K$ if $(i,k)=(r,k)=(r,i)=1$ and $i+k \equiv 0 \mod r$.
Then $i$ divides $t=(r-1)k-i$ only if $i=1$, so for $i>1$, these binomials are not
APN on an infinite number of extensions of $K$ by Corollary \ref{weilagain}.

\vspace{3mm}

\noindent {\bf Example 2.}
We now show that the sporadic
quadratic binomial function of \cite{EKP}:
    $$f = x^3 + u x^{36} \in GF(2^{10}),$$
    where $u \in S=\omega GF(2^5)^*\cup \omega^2 GF(2^5)^* $ and 
$\omega$ has order 3 in $GF(2^{10})$, is APN over at most a finite number of 
extensions of $GF(2^{10})$. 
If 
$$H_u(x,y) = 1+u\frac{x^7+y^7}{x+y}
(x^8 y + x y^8)^3$$ 
then $H_u(ax,ay) = H_{a^{33}u}(x,y)$, and $a^{33}u\in S$ if 
$u\in S$, so it suffices to prove
that $H_u(x,y)$ is absolutely irreducible for one $u\in S$.
Here $i=1$, $s=2$, and $t=3$ in our notation.
    
Let $\alpha$ be a root of the primitive polynomial $x^{10}+x^3+1$,
    and let $u=\alpha^{374}$.
    It can be easily checked (using a computer) that 
    $H_u(x,x^2+\alpha^{17}x+1)$ is irreducible in $K[x]$.
This implies that $H_u(x,y)$ is irreducible in $K[x,y]$.

The polynomial $H_u(x,y)$ has degree 33.
By Lemma \ref{rdivn}, if $H_u(x,y)$ is not absolutely irreducible
then it factors in one of the following three ways:
\medskip

-- 3 absolutely irreducible factors of degree 11 over $GF(2^{30})$ 

-- 11 absolutely irreducible factors of degree 3 over $GF(2^{110})$ 

-- 33 absolutely irreducible factors of degree 1 over $GF(2^{330})$.

\medskip
It is straightforward to check that $(1,0)$ is a point of 
multiplicity 3 on $H_u(x,y)$.
Since 11 and 33 do not divide 3, the second and third cases are not possible
by Lemma \ref{mults} applied with $r=1$, $m_0=3$.

Suppose finally that $H_u(x,y)$ factors over $GF(2^{30})$
into 3 factors of degree 11.
Again letting $\alpha$ be a root of the primitive polynomial $x^{10}+x^3+1$,
and $u=\alpha^{374}$, 
we checked with a computer that the polynomial 
$H_u(x,x^2+\alpha^{5}x+1)\in K[x]$ of degree $63$ has an irreducible
factor of degree $53$ in $K[x]$.  
Since 3 is relatively prime to 53, this factor remains irreducible over 
$GF(2^{30})$, and this is not compatible with 
the assumed factorization of $H_u(x,y)$.

\end{document}